\documentclass[a4paper,11pt]{amsart}
\usepackage[T1]{fontenc}
\usepackage{amsmath,amsthm,amscd,amsfonts,amsxtra,amssymb}

\title{A brief note on the second fundamental form of the variation of hodge structure of a smooth projective hypersurface}
\author{Emmanuel ALLAUD}
\email{eallaud@gmail.com}
\theoremstyle{plain}

\theoremstyle{definition}

\theoremstyle{remark}
\newtheorem*{rem*}{Remark}

\newcommand{\VHS}{variation of Hodge structure}

\newcommand{\IVHS}{infinitesimal variation of Hodge structure}

\newcommand{\TG}{Griffiths's transversality}

\newcommand{\C}{\mathbb C}

\newcommand{\der}[2]{\frac{\partial #1}{\partial #2}}

\newcommand{\FIXME}[1]{???}

\DeclareMathOperator{\End}{End}
\begin{document}
\begin{abstract}
We look here into the second order of a \VHS\ of smooth projective hypersurfaces. More precisely we give a simple algebraic way to describe the second fundamental form which in turn shows that the \VHS\ has to satisfy new second order partial differential equations.
\end{abstract}
\maketitle
We are interested in understanding the second order of a \VHS\ of the primitive cohomology associated to the universal family $\mathcal X$ of projective hypersurfaces:
$$
s\mapsto X_s=\{f_s=0\} \text{ with } f_s \in S^d(\C^{n+2})
$$
where $X=X_0$ is smooth.

A lot has already been done with respect to the \IVHS, see \cite{griffiths_periods} and \cite{voisin2002hodge}, \cite{voisin2003hodge} for example. Here we investigate the second fundamental form $II$ of the \VHS, and we get this description:
\begin{equation}
\forall u,v \in T_X \mathcal X,II(u,v)=u \circ v  
\end{equation}
where the operator is just usual composition of linear operators (\TG\ ensures that this is symmetric).

The main ingredient here is the identification of the \IVHS\ and its action on the cohomology with pieces of the jacobian ring of $X$ and its multiplication. We briefly recall here its definition and then prove our result.\\

Let $X=\{f=0\}$ be a smooth projective hypersurface of dimension $n$ and degree $d$, $f\in S^d(\C^{n+2})$ and the universal family of hypersurfaces $\mathcal X$:
$$
s\mapsto X_s=\{f_s=0\} \text{ with }f_s \in S^d(\C^{n+2}) \text{ and }f_0=f
$$
Let $S=S^d(\C^{n+2})$ and the jacobian ideal $J=\left( \der{f}{x_i}\right)$ then we define the Jacobian ring of $X$ by $$R=S/J.$$
The primitive cohomology of $X$ and the \IVHS\ $T=T_X \mathcal X$ are identified to some pieces of its jacobian ring via residues as described (for example) in \cite{voisin2003hodge}; this is summarized by the following commutative diagram:
\begin{equation} \label{diag_residus}
\begin{CD}
T @>>> \hom\left(H^{n-q,q},H^{n-q-1,q+1}\right) \\
@VVV   @VVV \\
R^d @>\times>> \hom\left(R^{(q+1)d-(n+2)},R^{(q+2)d-(n+2)}\right)
\end{CD}
\end{equation}
The vertical arrows are the residues isomorphisms (or the obvious maps induced by them), the upper arrow is the action of $T$ on $H^{n-p,p}$ and the lower one is the multiplication in $R$.

The identification is as follows (see \cite{donagi_torelli}): let $\alpha \in H^{n-p,p}(X)$  be mapped to $[P] \in R$ via residues, that is $\alpha$ is represented by the meromorphic form $\frac{P \Omega}{f^{p+1}}$, we have:
\begin{equation}\label{ivsh_action}
\der{\frac{P \Omega}{f^{p+1}_s}}{s}\mid_{s=0} = -(p+1)\frac{\der{f_s}{s}\mid_{s=0} P \Omega}{f^{p+2}}
\end{equation}
which is viewed as a class in $H^{n-(p+1),p+1}(X)$. This gives the usual identification of the \IVHS\ to the action of $R^d$ on $R$ by multiplication.

To establish the result about the second fundamental form of the \VHS\, we go on differentiating once more to check what happens at order 2; this time take $f_{s,t} \in S^d$:
\begin{align*}
\left( \frac{\partial ^2\frac{P \Omega}{f^{p+1}_{s,t}}}{\partial t \partial s} \right)_{s=0,t=0} &= -(p+1)\frac{\frac{\partial^2 f{s,t}}{\partial t \partial s}_{s=0,t=0} f^{p+2} - (p+2)\der{f_{s,t}}{s}_{s=0}\der{f_{s,t}}{t}_{t=0} f^{p+1}}{f^{2p+4}}  P \Omega\\
&=-(p+1)\left(\frac{\frac{\partial^2 f_{s,t}}{\partial t \partial s}_{s=0,t=0}P\Omega}{f^{p+2}}-(p+2) \frac{\der{f_{s,t}}{s}_{s=0}\der{f_{s,t}}{t}_{t=0}P\Omega}{f^{p+3}}\right)
\end{align*}
The first term is in the image of the \IVHS\ and then is quotiented out as we are talking about the 2nd fundamental form, the second term is in $H^{n-(p+2),p+2}(X)$ and is represented in the jacobian ring by the multiplication in $R$ by $\der{f_{s,t}}{s}\mid_{s=0}\der{f_{s,t}}{t}\mid_{t=0}$, that is $II\left(\der{f_{s,t}}{s}\mid_{s=0},\der{f_{s,t}}{t}\mid_{t=0} \right)$ is given by the composition of the action of $\der{f_{s,t}}{s}\mid_{s=0}$ and $\der{f_{s,t}}{t}\mid_{t=0}$ on $H^n(X)$.

Now to conclude we recall that the differential of the \VHS\ is given by a map (here $H=H(X)$):
$$
  T \to T_H D \cong \mathfrak g^-
$$
Where $\mathfrak g^-:=\bigoplus_{\substack{p=1}}^n \mathfrak g^{-p,p} $ is the Lie subalgebra of $\End(H,H)$ such that
$$\forall u \in \mathfrak g^{-p,p}, u(H^{n-k,k}) \subset H^{n-k-p,k+p}.$$
In fact $T$ is mapped into $\mathfrak g^{-1,1}$. The second fundamental form of the \VHS\ is then a map (we still denote by $T$ the image of $T$ in $T_HD$):
$$
II:S^2 T \to T_HD/T
$$
And what we just have seen translates in this setting to (up to a multiplicative constant):
$$
\forall u,v \in T_X \mathcal X,II(u,v)=u \circ v
$$

\begin{coronn}
This last equation means that the second fundamental form actually maps $S^2 T$ into $\mathfrak g^{-2,2}$ hence that its component in $\mathfrak g^{-1,1}/T \subset T_HD/T$ is zero, which in turn means that the \VHS\ satisfies a set of second order partial differential equations distinct from the ones coming from Griffiths' transversality.
\end{coronn}

\bibliographystyle{amsalpha}
\bibliography{article-1-en}
\end{document}